\documentclass[11pt]{amsart}
\numberwithin{equation}{section}
\newtheorem{theorem}{Theorem}

\newtheorem{corollary}{Corollary}

\numberwithin{theorem}{section} \numberwithin{lemma}{section}
\numberwithin{proposition}{section} \numberwithin{equation}{section}

\def\al{\aligned}
\def\eal{\endaligned}

\begin{document}

\tracingpages 1
\title[Sobolev]{\bf A uniform Sobolev inequality under Ricci flow}
\author{ Qi S. Zhang}

\address{Department of Mathematics,  University of California,
Riverside, CA 92521, USA }
\date{August 2007}

\maketitle







\tableofcontents


\begin{abstract}

Let ${\bf M}$ be a compact Riemannian manifold and the metrics
$g=g(t)$ evolve by the Ricci flow. We prove the following result.
The Sobolev imbedding by Aubin or Hebey, perturbed by a scalar
curvature term and modulo sharpness of constants, holds uniformly
for $({\bf M}, g(t))$ for all time if the Ricci flow exists for all
time; and if the Ricci flow develops a singularity in finite time,
then the same Sobolev imbedding holds uniformly after a standard
normalization. As a consequence,  long time non-collapsing results
are derived, which improve Perelman's local non-collapsing results.
Applications to Ricci flow with surgery are also presented.
\end{abstract}


\section{Introduction}


Let ${\bf M}$ be a compact Riemannian manifold of dimension $n \ge
3$ and $g$ be the metric. It is well known that a Sobolev
inequality of the following form holds: there exist positive
constants $A, B$
 such that, for all $v
\in W^{1, 2}({\bf M}, g)$,
\begin{equation}
\label{sob}
 \bigg( \int v^{2n/(n-2)} d\mu(g) \bigg)^{(n-2)/n} \le
A \int |\nabla v |^2 d\mu(g) + B \int v^2 d\mu(g).
\end{equation}

This inequality was proven by Aubin \cite{Au:1}  for
$A=K^2(n)+\epsilon$ with $\epsilon>0$ and $B$ depending on bounds
on the injectivity radius, sectional curvatures and derivatives.
Here $K(n)$ is the best constant in the Sobolev imbedding for
${\bf R}^n$. Hebey \cite{H:1} showed that $B$ can be chosen to
depend only on $\epsilon$, the injectivity radius and the lower
bound of the Ricci curvature. Hebey and Vaugon [HV] proved that
one can even take $\epsilon=0$. However the constant $B$ will also
depend on the derivatives of the curvature tensor.

Now consider a Ricci flow $({\bf M}, g(t))$,  $\quad \partial_t g
= - 2 Ric$. Due to the obvious importance of Sobolev inequality in
its analysis, it is highly desirable to have a uniform control on
the constants $A$ and $B$. For instance Sesum and Tian \cite{ST:1}
showed a uniform Sobolev imbedding for certain K\"ahler Ricci flow
under an additional lower bound of the Ricci curvatures. Chang and
Lu \cite{CL:1} proved that the Yamabe constant has nonnegative
derivative at the initial time under Ricci flow and under an
additional technical assumption. Also in the paper \cite{HV:2},
Hebey and Vaugon studied the evolution of Yamabe metrics  for
short time and  proved pinching results.

 Unfortunately the controlling geometric
quantities for $B$ as stated above are not invariant under the
Ricci flow in general. Nevertheless, by using a generalized
version Perelman's W entropy, we are able to prove that the above
Sobolev inequality (except for the sharp constants), with a
perturbation term involving the scalar curvature, is valid
uniformly under the Ricci flow or least after a suitable
normalization. As a consequence, we establish long time
non-collapsing result which generalizes Perelman's short time
result.

Here is the main result of the paper.

\begin{theorem}
\label{main} Let ${\bf M}$ be a compact Riemannian manifold with
dimension $n \ge 3$ and the metrics $g=g(t)$ evolve by the Ricci
flow $\partial_t g = - 2 Ric$. Then the following conclusions are
true.

(a). Suppose the Ricci flow exists for all time $t \in (0, \infty
)$. Then there exist positive functions $A, B$ depending only on the
initial metric $g(0)$ and $t$  such that, for all $v \in W^{1, 2}({\bf M},
g(t))$, $t>0$, it holds
\[
\bigg( \int v^{2n/(n-2)} d\mu(g(t)) \bigg)^{(n-2)/n} \le A \int (
|\nabla v |^2 + \frac{1}{4} R v^2 ) d\mu(g(t)) + B  \int v^2
d\mu(g(t)).
\]Here $R$ is the scalar curvature with respect to $g(t)$. The same
holds when $\infty$ is replaced by any $T>0$.

(b). Suppose the Ricci flow is smooth for $t \in (0, 1)$ and is
singular at $t=1$. Let $\tilde t = - \ln (1-t)$ and
$\tilde{g}(\tilde t) =\frac{1}{1-t} g(t)$ which satisfy a
normalized Ricci flow
\[
\partial_{\tilde t}\tilde g = - 2 \tilde{Ric} + \tilde g.
\]

Then there exist positive constants $A, B$ depending only on the
initial metric $g(0)$ such that, for all $v \in W^{1, 2}({\bf M},
\tilde g(\tilde t))$, $\tilde t>0$, it holds
\[
\bigg( \int v^{2n/(n-2)} d\mu(\tilde g(\tilde t)) \bigg)^{(n-2)/n}
\le A \int ( |\tilde \nabla v |^2 + \frac{1}{4} \tilde R v^2 )
d\mu(\tilde g(\tilde t)) + B \int v^2 d\mu(\tilde g(\tilde t)).
\]Here $\tilde R$ is the scalar curvature with respect to $\tilde g(\tilde t)$.
\end{theorem}
\medskip

{\it Remark}. (1). By the work of Hebey \cite{H:1} and Brouttelande
\cite{Br}, one can show that the constants in the theorem depending
on $g(0)$ only through $n$, lower bound of the Ricci curvature and
injectivity radius. This is explained in Step 2, Case (a) below.

(2). After the first version of the paper was posted on in the
arXiv, two interesting related papers [Y] and [HS] have appeared in
the arXiv. In one of them [Y], Professor Ye also pointed out an
inaccuracy in the statement of the previous Theorem 1.1 (a) and  its
proof. However, one can see in section 2 below, this can be easily
corrected by the same method of using the generalized Perelman
entropy. Actually this part of the proof can be done more simply by
using just Perelman's original W-entropy. In the papers [Y] and
[Hs], the controlling constants of the Sobolev imbedding were
computed more precisely and a full log Sobolev inequality were
deduced.


In case of K\"ahler Ricci flow, the theorem takes a particularly
simple form due to Perelman's boundedness result on the scalar
curvature, as explained in [ST].

\begin{corollary}
Let $({\bf M}, g(t))$ be the normalized K\"ahler Ricci flow
\[
\partial_t g_{i\bar j} = g_{i\bar j} - R_{i\bar j}
\]where ${\bf M}$ is a compact, K\"ahler manifold with complex dimension
$n$ and positive first Chern class. Then there exist positive
constants $A, B$ depending only on the initial metric $g(0)$ such
that, for all $v \in W^{1, 2}({\bf M}, g(t))$, $t>0$, it holds
\[
( \int v^{4n/(2n-2)} d\mu(g(t)) )^{(2n-2)/(2n)} \le A \int |\nabla
v |^2  d\mu(g(t)) + B \int v^2 d\mu(g(t)).
\]
\end{corollary}

By virtue of the work of Carron \cite{Ca} (see also Lemma 2.2 in
\cite{H:2} and its proof), we know that Sobolev imbedding in the
above form implies the following long term non-collapsing result for
geodesic balls on manifolds.

\begin{corollary}
In either case (a) or (b) in Theorem 1.1,  suppose also that the
scalar curvature is uniformly bounded. Then there exist two positive
constants $v_1$ and $v_2$ depending on only on $g(0)$, the bound on
the scalar curvature and on $t$ in case (a) such that
\[
| B(x, r) |_{g(t)} \ge \min [ v_1 r^n, v_2]; \qquad | B(x, r)
|_{\tilde g(\tilde t)} \ge \min [ v_1 r^n, v_2]
\]for all $ t \in [0, \infty)$ and $r>0$; and respectively
for all $ \tilde t \in [0, \infty)$ and $r>0$.
\end{corollary}

Using the property that the distances in a parabolic cube is
comparable if the curvature is bounded by the inverse radius square,
 this result generalizes Perelman's non-collapsing result from
short time to long time case.  Applications on local
non-collapsing with surgery are also possible. We will present it
in the addendum (last section) of the paper.

\section{Proof of Theorem 1.1}

{\bf Step 1.} We show that the mononotone property of $W$, the
 Perelman $W$ entropy,  implies  restricted Log Sobolev
inequalities (\ref{reslogsob1}) and (\ref{reslogsob2}) below.

{\it Case (a). We assume the Ricci flow exists for all time.}

 Clearly, by scaling in space time, we only need to prove the
result when $t = 1$.  Here is why. For a fixed  $t_0 \ge 1$. One
can do the scaling $t'=t/t_0$ and $g'=g/t_0$.

Let $t=1$. For any  $\epsilon \in (0, 1)$, we take
$\tau(\cdot)=\epsilon^2+1-\cdot$ so that $\tau_1 = 1+\epsilon^2$ and
$\tau_2=\epsilon^2$ (by taking $t_1=0$ and $t_2=1$).

Recall that Perelman's $W$ entropy is
\[
W(g, f, \tau)=\int_{\bf M}
                   \left(
                    \tau (  R+|\nabla f|^2)+f-n
                   \right)u\,d\mu(g(t))
\]where $u=\frac{e^{-f}}{(4 \pi \tau)^{n/2}}$.
Let $u_2$ be a minimizer of the entropy $W(g, f, \tau_2)$ for all
$u$ such that $\int u d\mu(g(t_2))=1$. We solve the backward heat
equation with the final value chosen as $u_2$ at $t=t_2$. Let $u_1$
be the value of the solution of the backward heat equation at
$t=t_1$. As usual, we define functions $f_i$ with $i=1, 2$ by the
relation $u_i = e^{-f_i}/(4 \pi \tau_i)^{n/2}$, $i=1, 2$. Then, by
the monotonicity of the $W$ entropy
\begin{equation}
\label{w1w2}
 \al inf_{\int u_0 d\mu(g(t_1))=1} W(g(t_1), f_0,
\tau_1) &\le
W(g(t_1), f_1, \tau_1) \le W(g(t_2), f_2, \tau_2) \\
&= inf_{\int u d\mu(g(t_2))=1} W(g(t_2), f, \tau_2).\eal
\end{equation}
Here $f_0$ and $f$ are given by the formulas
\[
u_0 = e^{-f_0}/(4 \pi \tau_1)^{n/2}, \qquad u = e^{-f}/(4 \pi
\tau_2)^{n/2}.
\]
Using these notations we can rewrite (\ref{w1w2}) as
\[
\al
 \inf_{\Vert u \Vert_1=1} \int_{\bf M}        &          \left(
                   \epsilon^2 (R+|\nabla \ln u |^2)- \ln u
                   -\ln(4\pi \epsilon^2)^{n/2}
                   \right) u\,d\mu(g(t))\\
&\ge \inf_{\Vert u_0 \Vert_1=1} \int_{\bf M} \left(
                   (1+\epsilon^2) (R+|\nabla \ln u_0 |^2)- \ln
                   u_0
                   -\ln(4\pi (1+\epsilon^2) )^{n/2}
                   \right) u_0\,d\mu(g(0)).
\eal
\]Observe that the $\ln(4\pi )^{n/2}$ terms on both sides of the
above inequality can be canceled. Denote $v=\sqrt{u}$ and
$v_0=\sqrt{u_0}$. We obtain, since $\epsilon \le 1$,
\begin{equation}
\label{vv0}
 \al
 \inf_{\Vert v \Vert_2=1} \int_{\bf M}         &          \left(
                   \epsilon^2 (R v^2 +4 |\nabla  v |^2)- v^2 \ln v^2
                   \right) \,d\mu(g(t)) - n \ln \epsilon\\
&\ge \inf_{\Vert v_0 \Vert_2=1} \int_{\bf M} \left(
                   (1+\epsilon^2) (R v^2_0+ 4|\nabla v_0 |^2)- v^2_0 \ln
                   v^2_0
                     \right) \,d\mu(g(0))-\ln 2^{n/2}.
\eal
\end{equation}

Since $({\bf M}, g(0))$ is a compact Riemannian
manifold, it is known that the following log Sobolev inequality
holds: {\it

Given $\epsilon \in (0, \sqrt{2 \pi}]$, for any $v_0 \in W^{1,
2}({\bf M})$ with $\Vert v_0 \Vert_2 =1$, there exists a positive
constant $L_0$ depending only on $g(0)$ such that
\begin{equation}
\label{logso}
 \int_{\bf M} v^2_0 \ln v^2_0 \,d\mu(g(0)) \le
\frac{\epsilon^2}{\pi} \int_{\bf M} |\nabla v_0 |^2 \,d\mu(g(0)) -
n \ln \epsilon + L_0.
\end{equation}
}Since the parameter $\epsilon$ is bounded from above by $\sqrt{2
\pi}$ while in the standard log Sobolev inequality $\epsilon$ is
unrestricted, we refer the above as restricted log Sobolev
inequality.

Actually even a sharp version of the log Sobolev inequality is valid
by [Br] Corollary 1.2. i.e.

{\it There exist positive constants $A_0$ and $B_0$ depending only
on the constants of the Sobolev imbedding (\ref{sob}) such that
for all $v_0 \in C^\infty({\bf M})$ verifying $\Vert v_0 \Vert_2 =
1$,
\[
\int_{\bf M} v^2_0 \ln v^2_0 \,d\mu(g(0)) \le \frac{1}{2} n \ln
\bigg( A_0 \int_{\bf M} |\nabla v_0 |^2 \,d\mu(g(0)) + B_0 \bigg).
\]}

This inequality is obtained by combining the Sobolev imbedding
(\ref{sob}) with H\"older's inequality with suitable powers and
taking the limit. If one does not care about the precise value of
$B_0$, one can even take $A_0=\frac{2}{n\pi e}$ (c.f. [Br]). For any
positive number $\xi_1$ and $\xi_2$ and $\delta \in (0, 1]$, it is
easy to see, that there exists $c>0$ such that,
\[
\ln(1+\xi_1+\xi_2) \le \delta (\xi_1+\xi_2) + |\ln \delta| +c \le
\delta \xi_1 + |\ln \delta| + \xi_2 + c.
\]
By choosing $\xi_1 = \int_{\bf M} |\nabla v_0 |^2 \,d\mu(g(0))$,
$\xi_2=B_0/A_0$ and taking $\delta$ appropriately, we know that
(\ref{logso}) is a direct consequence of this sharp log Sobolev
inequality.

 As a
consequence of this and the work \cite{H:1}, we know that $L_0$ can
be chosen to depend on $n$, the injectivity radius and the lower
bound of the Ricci curvature of $({\bf M}, g(0))$ only.

Substituting the log Sobolev inequality (\ref{logso}) to (\ref{vv0}), we
deduce
\begin{equation}
\label{vv1}
 \al
 \inf_{\Vert v \Vert_2=1} \int_{\bf M}         &          \left(
                   \epsilon^2 (R v^2 +4 |\nabla  v |^2)
                   - v^2 \ln v^2
                   \right) \,d\mu(g(t)) -n \ln \epsilon \\
&\ge - \max R^-(\cdot, 0) - L_0
 -c.
\eal
\end{equation}Therefore, by renaming $\epsilon$ and $c$, we
reach the uniform restricted log Sobolev inequality:
\begin{equation}
\label{reslogsob1}
 \int_{\bf
M}   v^2 \ln v^2 \,d\mu(g(t)) \le \frac{\epsilon^2}{2\pi}
\int_{\bf M} \big( |\nabla  v |^2 + \frac{1}{4} R v^2 \big)
d\mu(g(t))
  - n \ln \epsilon + L + \max R^-(\cdot, 0).
\end{equation}
Here $L$ depends only on $n$ and $g(0)$ through the lower bound of
the Ricci curvature and injectivity radius.

\medskip

 {\it Case (b). We assume the Ricci flow exists for $t \in
(0, 1)$ and becomes singular at $t=1$.}

We have
\begin{equation}
\label{reslogsob2}
 \int_{\bf M} v^2 \ln v^2 \,d\mu(\tilde g(\tilde t))
\le \frac{\epsilon^2}{2\pi} \int_{\bf M} \big( |\tilde \nabla  v |^2
+ \frac{1}{4}\tilde R v^2 \big) d\mu(\tilde g(\tilde t))
  - n \ln \epsilon + L + \max R^-(\cdot, 0).
\end{equation}Here $L$ depends only on $n$ and $g(0)$ through the
lower bound of the Ricci curvature and injectivity radius.

This  is an immediate consequence of Case (a) by scaling.

\medskip

{\bf Step 2.}

Fix a time $t_0$ during the Ricci flow or  $\tilde t_0$ during
 the normalized one. Suppose on $({\bf M}, g(t_0))$ or $({\bf M},
\tilde g( \tilde t_0))$,  the restricted Log Sobolev inequalities
(\ref{reslogsob1}) and (\ref{reslogsob2}) hold respectively. We
show that they imply short time heat kernel upper bound for the
fundamental solution of
\begin{equation}
\label{fixeq}
 \Delta u(x, t) - \frac{1}{4} R(x, t_0) u(x, t) -\partial_t u(x, t)=0
\end{equation}under the fixed
metric $g(t_0)$ or $\tilde g(\tilde t_0)$. The proof, which does
not distinguish between the Ricci flow or the normalized Ricci
flow case, follows the original ideas of Davies \cite{Da}. There
are only two extra issues to deal with here. One is that the range
of $\epsilon$ is restricted. The other is that the negative part
of the scalar curvature may make the semigroup generated by
$\Delta -\frac{1}{4} R$ not contractive. However the modification
in the proof is moderate since we are dealing with short time
upper bound only. We also benefit from the fact that the most
negative value of the scalar curvature does not decrease under
either the Ricci flow or the normalized one in the theorem. For
this reason, we will be brief in the presentation.

Let $u$ be a positive solution to (\ref{fixeq}). Given $T \in (0,
1]$ and $t \in (0, T)$, we take $p(t)=T/(T-t)$ so that $p(0)=1$
and $p(T)=\infty$. By direct computation
\[
\al
\partial_t \Vert u \Vert_{p(t)} &= \partial_t \bigg( \int_{\bf M}
u^{p(t)}(x, t) dx \bigg)^{1/p(t)}\\
&=-\frac{p'(t)}{p^2(t)} \Vert u \Vert_{p(t)} \ln \int_{\bf M}
u^{p(t)}(x, t) dx
    +\frac{1}{p(t)} \bigg( \int_{\bf M} u^{p(t)}(x, t) dx \bigg)^{(1/p(t))-1}
          \\
          &\qquad \times \bigg[ \int_{\bf M} u^{p(t)} (\ln u) p'(t) dx +
    p(t) \int_{\bf M} u^{p(t)-1} ( \Delta u - \frac{1}{4} R u ) dx \bigg].
\eal
\]Here $dx$ means the integral element with respect to $g(t_0)$.
We adopt this notation to symbolize the property that $g(t_0)$ is
not evolving with respect to $t$. Using integration by parts on
the term containing $\Delta u$ and multiplying both sides by
$p^2(t) \Vert u \Vert^{p(t)}_{p(t)}$, we reach
\[
\al &p^2(t)  \Vert u \Vert^{p(t)}_{p(t)}
\partial_t \Vert u \Vert_{p(t)}\\
&= -p'(t) \Vert u \Vert^{p(t)+1}_{p(t)} \ln \int_{\bf M}
u^{p(t)}(x, t) dx
  + p(t) \Vert u \Vert_{p(t)} p'(t) \int_{\bf M}
      u^{p(t)} \ln u (x, t) dx\\
      &\qquad
         -p^2(t)(p(t)-1) \Vert u \Vert_{p(t)} \int_{\bf M}
            u^{p(t)-2} |\nabla u|^2(x, t) dx\\
            &\qquad \qquad
              - p^2(t) \Vert u \Vert_{p(t)} \int_{\bf M} \frac{1}{4} R(x, t_0)
            u^{p(t)}(x, t) dx.
\eal
\]Dividing both sides by $ \Vert u \Vert_{p(t)}$, we obtain
\[
\al
&
p^2(t)  \Vert u \Vert^{p(t)}_{p(t)}
\partial_t \ln \Vert u \Vert_{p(t)}\\
&= -p'(t) \Vert u \Vert^{p(t)}_{p(t)} \ln \int_{\bf M} u^{p(t)}(x,
t) dx
  + p(t)  p'(t) \int_{\bf M}
      u^{p(t)} \ln u (x, t) dx\\
      &\qquad
         -4(p(t)-1)  \int_{\bf M}
            |\nabla (u^{p(t)/2}) |^2(x, t) dx
              - p^2(t)  \int_{\bf M} \frac{1}{4} R(x, t_0)
            (u^{p(t)/2})^2(x, t) dx.
\eal
\]Merging the first two terms on the righthand side of the above
equality and making the substitution $v = u^{p(t)/2}/\Vert
u^{p(t)/2} \Vert_2$, we arrive at, after dividing by $\Vert u
\Vert^{p(t)}_{p(t)}$,
\[
\al & p^2(t)
\partial_t \ln \Vert u \Vert_{p(t)}\\
&= p'(t)   \int_{\bf M} v^2 \ln v^2 (x, t) dx
         -4(p(t)-1)  \int_{\bf M}
            |\nabla v |^2(x, t) dx
              - p^2(t)  \int_{\bf M} \frac{1}{4} R(x, t_0) v^2(x, t) dx\\
      &= p'(t)   \int_{\bf M} v^2 \ln v^2 (x, t) dx
         -4(p(t)-1)  \int_{\bf M}
            (|\nabla v |^2(x, t) +\frac{1}{4} R(x, t_0) v^2) dx\\
            &\qquad \qquad
              + (4(p(t)-1)- p^2(t))  \int_{\bf M} \frac{1}{4} R(x, t_0) v^2(x, t)
              dx.
 \eal
\]It is easy to check  $\Vert v \Vert_2 = 1$ and also
\[
\frac{4(p(t)-1)}{p'(t)} = \frac{4 t (T-t)}{T} \le T \le 1,
\]
\[
-T \le \frac{4(p(t)-1)-p^2(t)}{p'(t)} = \frac{4t(T-t)-T^2}{T} \le 0.
\]Hence
\[
\al & p^2(t)
\partial_t \ln \Vert u \Vert_{p(t)} \\
&\le p'(t) \bigg(  \int_{\bf M} v^2 \ln v^2 (x, t) dx -
\frac{4(p-1)}{p'(t)} \int_{\bf M}
            (|\nabla v |^2(x, t) +\frac{1}{4} R(x, t_0) v^2) dx
+ T \sup R^-(x, t_0) \bigg). \eal
\]Take $\epsilon$ so that
\[
\frac{\epsilon^2}{\pi} = \frac{4(p(t)-1)}{p'(t)} \le T \le 1 \] in
the restricted log Sobolev inequality (\ref{reslogsob1}), we
deduce
\[
p^2(t)
\partial_t \ln \Vert u \Vert_{p(t)} \le
 p'(t) \bigg( -n \ln \sqrt{\pi 4(p(t)-1)/p'(t)} + L + T \sup
 R^-(x, 0) \bigg).
 \]Here we also used the fact that $\sup
 R^-(x, t_0) \le \sup
 R^-(x, 0)$ as remarked earlier.

 Note that $p'(t)/p^2(t)=1/T$ and $4(p(t)-1)/p'(t)=4t (T-t)/T$.
 Hence
 \[
\partial_t \ln \Vert u \Vert_{p(t)} \le
 \frac{1}{T} \bigg( - \frac{n}{2} \ln \pi 4t (T-t)/T + L + T \sup
 R^-(x, 0) \bigg).
 \]This yields, after integration from $t=0$ to $t=T$,
 \[
 \ln  \frac{\Vert u(\cdot, T) \Vert_\infty}{\Vert u(\cdot, 0)
 \Vert_1}
 \le - \frac{n}{2} \ln (4 \pi T) + L + T \sup
 R^-(x, 0).
 \]Since
 \[
 u(x, T) = \int_{\bf M}  P(x, y, T) u(y, 0) dy,
 \]
this shows
 \[
 P(x, y, T) \le \frac{\exp(L+ T \sup
 R^-(x, 0))}{(4 \pi T)^{n/2}}.
 \]Here $P$ is heat kernel of (\ref{fixeq}).
\medskip

{\bf Step 3.} We show that the short term heat kernel upper bound
implies the Sobolev imbedding in Theorem 1.1.

\medskip
 This is more or less
standard. Let $t_0$ be a fixed time during Ricci flow. Let $F=\sup
 R^-(x, 0)$ and $P_F$ be the heat kernel of the operator
 $\Delta - \frac{1}{4} R(x, t_0) -F-1$. Since $R^-(x, t_0) \le F$, from the short time upper
 bound for $P$, we know that
 $P_F$ obeys the global upper bound
 \[
 P_F(x, t, y) \le \frac{\Lambda}{t^{n/2}}, \qquad t >0.
 \]Here $\Lambda$ depends only on $L$ and $F$.
Moreover $P_F$ is a contraction. By H\"older inequality, for any $f
\in L^2({\bf M})$, we have
\[
|\int_{\bf M} P_F(x, t, y) f(y) dy| \le \bigg( \int_{\bf M} P^2_F(x,
t, y) dy \bigg)^{1/2} \ \Vert f \Vert_2 \le \Lambda^{1/2} t^{-n/4}
\Vert f \Vert_2.
\]
The Sobolev inequality in Theorem 1.1 now follows from Theorem
2.4.2 in \cite{Da}. i.e. there exist positive constants $A, B$
depending only on the initial metric through $\Lambda$ such that,
for all $v \in W^{1, 2}({\bf M}, g(t_0))$, it holds
\[
\bigg( \int v^{2n/(n-2)} d\mu(g(t_0)) \bigg)^{(n-2)/n} \le A \int (
|\nabla v |^2 + \frac{1}{4} R v^2 ) d\mu(g(t_0)) + B \int v^2
d\mu(g(t_0)).
\]The same also holds for the normalized Ricci flow. Since $t_0$
is arbitrary, the proof is done.

 One can also prove it by establishing a Nash
type inequality first and using an argument in \cite{BCLS}. \qed

 \emph{Acknowledgement.}  We would like to thank Professors T.
Coulhon, L. F. Tam, G. Tian and X. P. Zhu for their helpful
communications.
\vfill
\eject

\section{ADDENDUM}

In this section, we generalize the previous result to Ricci flow
with surgery. As an application we prove a noncollapsing result with
surgery without using the concepts of reduced distance and reduced
volume by Perelman.

Here is the main result. For detailed information and related
terminology such as $(r, \delta)$ surgery, $\delta$ neck etc on
Ricci flow with surgery we refer the reader to [CZ], [KL] and [MT].
\medskip

{\it {\bf Theorem A.1.} Given real numbers $T_1<T_2<T_3$, let $({\bf
M}, g(t))$ be a 3 dimensional Ricci flow with normalized initial
condition defined on the time interval containing $[T_1, T_3]$.
Suppose the following conditions are met.

(a). $g(t)$ is smooth except when $t=T_2$.

(b). A $(r, \delta)$ surgery occurs at $t=T_2$ with parameter
$h$. i.e. the surgery occurs in a $\delta$ neck of radius $h$. Here
$\delta$ is sufficiently small.

(c). For $t \in [T_1, T_2)$ and $A_1>0$, the Sobolev imbedding
\[
\bigg( \int v^{2n/(n-2)} d\mu(g(t)) \bigg)^{(n-2)/n} \le A_1 \int (
|\nabla v |^2 + \frac{1}{4} R v^2 ) d\mu(g(t)) + A_1  \int v^2
d\mu(g(t)). \
\]holds for all $v \in W^{1, 2}({\bf M}, g(t))$. Here $R$ is the
scalar curvature.

Then for all $t \in (T_2, T_3]$, the Sobolev imbedding below holds
for all $v \in W^{1, 2}({\bf M}, g(t))$.
\[
\bigg( \int v^{2n/(n-2)} d\mu(g(t)) \bigg)^{(n-2)/n} \le A_2 \int (
|\nabla v |^2 + \frac{1}{4} R v^2 ) d\mu(g(t)) + A_2  \int v^2
d\mu(g(t))
\]
Here
\[
A_2 \le a_1 e^{a_2(T_3-T_2 + R^-_0)} ( A_1 +  R^-_0 + a_3)
\]with $a_i, i=1, ..., 3$ being positive numerical constants; and
$R^-_0=\sup R^-(x, 0).$

Moreover, the Ricci flow is $\kappa$ noncollapsed in the whole
interval $[T_1, T_3]$  under a scale $\epsilon$ where $\kappa$ and
$\epsilon$ are independent of the surgery radius $h$, or the
canonical neighborhood scale $r(t)$ in $(T_2, T_3]$ or the smallness
of $\delta$. }
\medskip

{\it Remark.}  1. It is well known that any compact Riemannian manifold
with a given metric supports a Sobolev imbedding ([Au], [H1] e.g).
Therefore the theorem shows that under a Ricci flow with finite
number of surgeries in finite time, the Sobolev imbedding is
preserved. By Lemma A.2 below, $\kappa$ noncollapsing property is
also preserved in this case. The constant $\kappa$ will depend on
time and number of surgeries, but it is independent of the surgery
scale or the canonical neighborhood scale $r(t)$.

2.  In this paper we use the following definition of $\kappa$ non-collapsing by Perelman [P2],  as elucidated in Definition 77.9 of [KL].

Let ${\bf M}$ be a Ricci flow with surgery defined on $[a, b]$. Suppose that $(x_0, t_0) \in {\bf M}$ and $r>0$ are such that $t_0-r^2 \ge a$, $B(x_0, t_0, r)
\subset {\bf M}^-_{t_0}$ is a proper ball and the parabolic ball $P(x_0, t_0, r, -r^2)$
is unscathed. Then ${\bf M}$ is $\kappa$-collapsed at $(x_0, t_0)$ at scale $r$
if $|Rm| \le r^{-2}$ on $P(x_0, t_0, r, -r^2)$ and $vol(B(x_0, t_0, r)) < \kappa r^3$;
otherwise it is $\kappa$-noncollapsed.

\medskip
In order to prove the theorem, we need the following two lemmas. The
first lemma extends Theorem 1.1 in that  a more accurate upper bound
for the evolving Sobolev constants are derived.
\medskip

{\it {\bf Lemma A1.}  Let $({\bf M}, g(t))$ be a smooth Ricci flow for $t \in [0, T_0]$.
Suppose for a constant $A(0)>0$, the following Sobolev imbedding holds for the initial metric:

for all $v \in W^{1, 2}({\bf M},
g(0))$,
\[
\bigg( \int v^{2n/(n-2)} d\mu(g(0)) \bigg)^{(n-2)/n} \le A (0)\int (
|\nabla v |^2 + \frac{1}{4} R v^2 ) d\mu(g(0)) + A(0)  \int v^2
d\mu(g(0)).
\leqno(A.1)
\]

Then there exists a positive constant $A(T_0)$ such that
the following Sobolev imbedding holds for the  metric $g(T_0)$:

for all $v \in W^{1, 2}({\bf M},
g(T_0))$,
\[
\bigg( \int v^{2n/(n-2)} d\mu(g(T_0)) \bigg)^{(n-2)/n} \le A (T_0)\int (
|\nabla v |^2 + \frac{1}{4} R v^2 ) d\mu(g(T_0)) + A(T_0)  \int v^2
d\mu(g(T_0)).
\leqno(A.2)
\]
Moreover, there exists numerical constants $a_1, a_2>0$ such that
\[
A(T_0) \le a_1 \exp (a_2 (T_0 +1 + R^-_0)) A(0),
\]where $R^-_0 = \sup R^-(x, 0)$.

In the above and later, $R$ stands for the scalar curvature under the given metrics.

}

\noindent {\bf Proof.}

We begin by deriving a log Sobolev inequality from the Sobolev inequality (A.1) using the
argument in [Br], p119-120.

Using (A.1) and interpolation (H\"older) inequality, we have (see
(1.3) in [Br])

\[
\al \left( \int v^r d\mu(g(0)) \right)^{2/(r\theta)} \le & \left( A
(0)\int ( |\nabla v |^2 + \frac{1}{4} R v^2 ) d\mu(g(0)) +A(0)  \int
v^2
d\mu(g(0)) \right) \\
&\qquad \times \left(  \int v^s d\mu(g(0)) \right) ^{2(1-\theta)/(s \theta)}
\eal
\]where
\[
\frac{1}{r} = \frac{\theta (n-2)}{2n} + \frac{1-\theta}{s}.
\]

Next we take $r=2$ and $v$ such that $\Vert v \Vert_2 =1$, divide the last integral term on the right hand side to the left,
and take log on the resulting inequality. Then we take $\theta \to 0$, i.e. $s \to 2$.
It is straight forward to check that we will arrive at the subsequent log Sobolev inequality:

 for those $v \in W^{1, 2}({\bf M}, g(0))$ such that $\Vert v \Vert_2 =1$, it holds
\[
\int v^2 \ln v^2 d\mu (g(0)) \le \frac{n}{2} \ln
\left( A (0)\int (
|\nabla v |^2 + \frac{1}{4} R v^2 ) d\mu(g(0)) + A(0)   \right).
\leqno(A.3)
\]

Recall the elementary inequality: for all $z, \sigma >0$,
\[
\ln z \le \sigma z - \ln \sigma -1.
\leqno(A.4)
\]Its proof goes as follows.  Write
\[
f(z)=\ln z -\sigma z + \ln \sigma +1.
\]Then $f'(z)=(1/z)-\sigma$ and $f^{''}(z)=-1/z^2<0$. So $f$ is concave down. It is absolute
maximum is reached when $f'=0$, i.e. when $z=1/\sigma$. But
$
f(1/\sigma) =0
$. Hence  $f(z) \le 0$ for all $z>0$, yielding (A.4).

Using (A.4) on the right hand side of (A.3), we deduce
\[
\int v^2 \ln v^2 d\mu (g(0)) \le \frac{n}{2} \sigma
\left( A (0)\int (
|\nabla v |^2 + \frac{1}{4} R v^2 ) d\mu(g(0)) + A(0)   \right) - \frac{n}{2} \ln \sigma -
\frac{n}{2}.
\leqno(A.5)
\]We write
\[
\frac{n}{2} \sigma A(0)  \equiv \frac{\alpha^2}{2 \pi}, \quad i.e. \quad
\sigma = \frac{\alpha^2}{2 \pi} \frac{2}{n A(0)}.
\]Then (A.5) becomes
\[
\int v^2 \ln v^2 d\mu (g(0)) \le \frac{\alpha^2}{2 \pi} \int (
|\nabla v |^2 + \frac{1}{4} R v^2 ) d\mu(g(0)) +  \frac{\alpha^2}{2 \pi}
- \frac{n}{2} \ln  \frac{\alpha^2}{2 \pi} + \frac{n}{2} \ln \frac{n A(0)}{2} -
\frac{n}{2}.
\leqno(A.6)
\]

Our next task is to use the monotonicity of Perelman's W entropy to extend (A.6) to a
log Sobolev inequality for $g(T_0)$.

Let $t \in (0, T_0]$ and $\epsilon \in (0, 1]$. We define
\[
\tau = \tau(t) = \epsilon^2 + T_0 -t
\]
so that $\tau_1 = \epsilon^2+T_0$ and
$\tau_2=\epsilon^2$ (by taking $t_1=0$ and $t_2=T_0$).

Let $u_2$ be a minimizer of the entropy $W(g, f, \tau_2)$ for all
$u$ such that $\int u d\mu(g(t_2))=1$. We solve the backward heat
equation with the final value chosen as $u_2$ at $t=t_2$. Let $u_1$
be the value of the solution of the backward heat equation at
$t=t_1$. As usual, we define functions $f_i$ with $i=1, 2$ by the
relation $u_i = e^{-f_i}/(4 \pi \tau_i)^{n/2}$, $i=1, 2$. Then, by
the monotonicity of the $W$ entropy
\[
 \al inf_{\int u_0 d\mu(g(t_1))=1} W(g(t_1), f_0,
\tau_1) &\le
W(g(t_1), f_1, \tau_1) \le W(g(t_2), f_2, \tau_2) \\
&= inf_{\int u d\mu(g(t_2))=1} W(g(t_2), f, \tau_2).\eal
\]
Here $f_0$ and $f$ are given by the formulas
\[
u_0 = e^{-f_0}/(4 \pi \tau_1)^{n/2}, \qquad u = e^{-f}/(4 \pi
\tau_2)^{n/2}.
\]
Using these notations we can rewrite the above as
\[
\al
 \inf_{\Vert u \Vert_1=1} \int_{\bf M}        &          \left(
                   \epsilon^2 (R+|\nabla \ln u |^2)- \ln u
                   -\ln(4\pi \epsilon^2)^{n/2}
                   \right) u\,d\mu(g(T_0))\\
&\ge \inf_{\Vert u_0 \Vert_1=1} \int_{\bf M} \left(
                   (\epsilon^2+T_0) (R+|\nabla \ln u_0 |^2)- \ln
                   u_0
                   -\ln(4\pi (\epsilon^2 +T_0))^{n/2}
                   \right) u_0\,d\mu(g(0)).
\eal
\]Denote $v=\sqrt{u}$ and $v_0 =\sqrt{u_0}$. This inequality is converted to
\[
\al
 &\inf_{\Vert v \Vert_2=1} \int_{\bf M}                \left(
                   \epsilon^2 (R v^2+4 |\nabla v |^2)- v^2 \ln v^2
                   \right) \,d\mu(g(T_0)) - \ln(4\pi \epsilon^2)^{n/2}\\
&\ge \inf_{\Vert v_0 \Vert_2=1} \int_{\bf M} \left(
                   4(\epsilon^2+T_0) (\frac{1}{4}R v^2_0+|\nabla v_0 |^2)- v^2_0 \ln
                   v^2_0
                   \right) \,d\mu(g(0)) -\ln(4\pi (\epsilon^2 +T_0))^{n/2}.
\eal
\]

Apply the log Sobolev inequality (A.6) on the right hand side of the
above inequality with the choice
\[
\frac{\alpha^2}{2 \pi} = 4 (T_0 + \epsilon^2).
\]We arrive at
\[
\al
 \inf_{\Vert v \Vert_2=1} \int_{\bf M}        &          \left(
                   \epsilon^2 (R v^2+4 |\nabla v |^2)- v^2 \ln v^2
                   \right) \,d\mu(g(T_0)) - n \ln \epsilon\\
&\ge - 4(\epsilon^2+T_0) - \frac{n}{2} \ln \frac{n A(0)}{2}-c.
\eal
\]Here $c$ is a numerical constant.
Since $\epsilon \le 1$, the above implies the restricted log Sobolev
inequality for $g(T_0)$:
\[
 \int_{\bf M} v^2 \ln v^2 \,d\mu(g(T_0)) \le
 \frac{\epsilon^2}{2 \pi} \int_{\bf M}
                   ( |\nabla v |^2 + \frac{1}{4} R v^2) \,d\mu(g(T_0))
                   -n \ln \epsilon + 4 (T_0+1) + \frac{n}{2} \ln
                   A(0)
                   + c.
 \leqno(A.7)
 \]

 Let $p(x, y, t)$ be the heat kernel of $\Delta - \frac{1}{4} R$ in
 $({\bf M}, g(T_0))$. Then (A.7) implies, for $t \in (0, 1]$,
 \[
 p(x, y, t) \le \exp ( 4(T_0+1) +\frac{n}{2} \ln
                   A(0) + c + R^-_0) \frac{1}{(4 \pi t)^{n/2}} \equiv
                   \frac{\Lambda}{t^{n/2}}.
 \leqno(A.8)
 \]This can be proven by adapting the standard method in heat
  kernel estimate as demonstrated in section 2 of the paper.
It is well known that (A.8) implies the desired Sobolev imbedding
for $g(T_0)$, i.e.

for all $v \in W^{1, 2}({\bf M}, g(T_0))$,
\[
\bigg( \int v^{2n/(n-2)} d\mu(g(T_0)) \bigg)^{(n-2)/n} \le A
(T_0)\int ( |\nabla v |^2 + \frac{1}{4} R v^2 ) d\mu(g(T_0)) +
A(T_0)  \int v^2 d\mu(g(T_0)).
\]
By keeping track of the constants (see also Theorem 2.2 in [S-C]),
we know that
\[
A(T_0) \le c_1 \Lambda^{2/n}
\]for some numerical constant $c_1$.
Therefore, there exists numerical constants $a_1, a_2>0$ such that
\[
A(T_0) \le a_1 \exp (a_2 (T_0 +1 + R^-_0)) A(0).
\]\qed

\medskip

The second Lemma of the addendum relates the Sobolev imbedding in
Lemma A.1 to local noncollapsing of volume of geodesic balls. We
follow the idea in [Ca].
\medskip

{\it {\bf Lemma A.2.} Let $({\bf M}, g)$ be a Riemannian manifold.
Given $x_0 \in {\bf M}$ and $r \in (0, 1]$. Let $B(x_0, r)$ be a
proper geodesic ball, i.e. ${\bf M} - B(x_0, r)$ is non empty.
Suppose the scalar curvature $R$ satisfies $|R(x)| \le 1/r^2$ in
$B(x_0, r)$ and the following Sobolev imbedding holds:
 for all $v \in W^{1, 2}_0(B(x_0, r))$, and a constant $A \ge 1$,
 \[
\bigg( \int v^{2n/(n-2)} d\mu(g) \bigg)^{(n-2)/n} \le A \int (
|\nabla v |^2 + \frac{1}{4} R v^2 ) d\mu(g) + A  \int v^2 d\mu(g).
\]Then $|B(x_0, r)| \ge 2^{-(n+5) n/2} A^{-n/2} r^n.$
}
\medskip

\noindent {\bf Proof.} Since $R \le 1/r^2$, $r \le 1$ and $A \ge 1$
by assumption, the Sobolev imbedding can be simplified to
\[
\bigg( \int v^{2n/(n-2)} d\mu(g) \bigg)^{(n-2)/n} \le A \int |\nabla
v |^2  d\mu(g) + \frac{2 A}{r^2}  \int v^2 d\mu(g).
\]

Under the scaled metric $g_1 = g/r^2$, we have, for all $v \in W^{1,
2}_0(B(x_0, 1, g_1))$,
\[
\bigg( \int v^{2n/(n-2)} d\mu(g_1) \bigg)^{(n-2)/n} \le A \int
|\nabla v |^2  d\mu(g_1) + 2 A  \int v^2 d\mu(g_1).
\]Now, by [C] (see p33, line 4 of [H2]), it holds
\[
|B(x_0, 1, g_1)|_{g_1} \ge \min \{\frac{1}{2 \sqrt{2A}},
\frac{1}{2^{(n+4)/2} \sqrt{2A}} \}^n.
\]Therefore
\[
|B(x_0, r,, g)|_g \ge 2^{-(n+5) n/2} A^{-n/2} r^n.
\]\qed
\medskip

Now we are in a position to give a

\medskip

{\bf Proof of Theorem A.1.}

We begin by finding an upper bound for the Sobolev constant of for
the metric right after the surgery.  Recall that surgery occurs at a
$\delta$ neck, called ${\bf N}$, of radius $h$ such that $({\bf N},
h^{-2} g)$ is $\delta$ close in the $C^{[\delta^{-1}]}$ topology to
the standard round neck $S^2 \times (-\delta^{-1}, \delta^{-1})$ of
scalar curvature $1$. Let $\Pi$ be the diffeomorphism from the
standard round neck to ${\bf N}$ in the definition of $\delta$
closeness. Denote by $z$ for a number in $(-\delta^{-1},
\delta^{-1})$.  For $\theta \in S^2$, $(\theta, z)$ is a
parametrization of ${\bf N}$ via the diffeomorphism $\Pi$. We can
identify the metric on ${\bf N}$ with its pull back on the round
neck by $\Pi$ in this manner.

Following  the notations on p424 of [CZ], the metric $\tilde
g=\tilde g(T_2)$ right after the surgery is given by
\[
\tilde g =
\begin{cases}
\bar g, \qquad z \le 0,\\
e^{-2 f} \bar g, \qquad z \in [0, 2],\\
\phi e^{-2 f} \bar g + (1-\phi)e^{-2 f} h^2 g_0 , \qquad z \in [2, 3],\\
e^{-2 f} h^2 g_0 , \qquad z \in [3, 4].
\end{cases}
\]Here $\bar g$ is the nonsingular part of the  $\lim_{t \to
T^-_2} g(t)$; $g_0$ is the standard metric on the round neck; and
$f=f(z)$ is a smooth function given by (c.f. p424 [CZ])
\[
\al
 &f(z)=0, \, z \le 0; \quad f(z) = ce^{-P/z}, \, z \in (0, 3]; \,
f''(z)>0, \, z \in [3, 3.9]; \\
&\quad f(z)=-\frac{1}{2} \ln (16-z^2), \, z \in [3.9, 4]. \eal
\]Here a small $c>0$ and a large $P>0$ are suitably chosen to
ensure that the Hamilton-Ivey pinching condition remains valid.
$\phi$ is a smooth bump function with $\phi=1$ for $z \le 2$ and
$\phi=0$ for $z \ge 3$.

Now we define another cut-off function $\eta$ on ${\bf N}$ such that, $0 \le \eta \le 1$
and
\[
\eta(x) =
\begin{cases} 1, \quad -100 \le z \le -10,\\
0, \quad z \ge 0.
\end{cases}
\]Here $x \in \Pi (S^2 \times \{ z \})$. i.e. $\eta(x) \equiv \xi(z(x))$ where
$\xi$ is a suitable cut-off function on the real line. We then
extend $\eta$ to be a Lipschitz cut-off function on the whole
manifold $({\bf M}, \tilde g(T_2))$ by taking $\eta(x)=1$ when $x
\in {\bf M} - \Pi(S^2 \times (-100, 4]).$

By $\delta$ closeness, there exists a constant $c$, independent of
the smallness of $\delta$, such that $|\nabla \eta (x)| \le c
h^{-1}.$    Here $\nabla$ is with respect to the metric $\tilde
g(T_2)$ which is identical to  the unscathed metric $\bar g(T_2)$ in
the support of $\nabla \eta$.

Notice that $\nabla \eta (x)=0$ when $\eta(x)=1$; and when $\eta(x)<1$, then $x$ is in the
$\delta$ neck where the scalar curvature is close to $h^{-2}$ when $\delta$ is small.
Therefore there exists a constant $c$, independent of the smallness of $\delta$, such that
\[
|\nabla \eta (x) |^2 \le c ( R(x, T_2) + R^-_0), \qquad x \in {\bf
M}, \, \text{with metric} \, \tilde g(T_2). \leqno(A.9)
\]Here, as before, $R^-_0 = \sup R^-(x, 0) \ge \sup R^-(x, t)$, $t \ge 0$.

Now we pick $v \in C^\infty({\bf M}, \tilde g(T_2))$. We compute
\[
\al
&\left( \int v^{2n/(n-2)} d\mu(\tilde g(T_2)) \right)^{(n-2)/n}   =
\left( \int [( \eta + 1-\eta)v]^{2n/(n-2)} d\mu(\tilde g(T_2)) \right)^{(n-2)/n}\\
&\le 2 \left( \int (\eta v)^{2n/(n-2)} d\mu(\tilde g(T_2)) \right)^{(n-2)/n} +
2 \left( \int ((1-\eta) v)^{2n/(n-2)} d\mu(\tilde g(T_2)) \right)^{(n-2)/n} \\
&\equiv 2 I_1 + 2 I_2.
\eal
\leqno(A.10)
\]

Observe that $supp (\eta v)$ is in the unscathed portion of ${\bf M}$ so that $\tilde g(T_2)=
\bar g(T_2)$ here. Hence condition (c) in the statement of the theorem is applicable to $\eta v$. Therefore
\[
\al
I_1 &\le A_1 \int (
|\nabla (\eta v) |^2 + \frac{1}{4} R (\eta v)^2 ) d\mu(\bar g(T_2)) + A_1  \int (\eta v)^2
d\mu(\bar g(T_2))\\
&=A_1 \int (
|\nabla (\eta v) |^2 + \frac{1}{4} R (\eta v)^2 ) d\mu(\tilde g(T_2)) + A_1  \int (\eta v)^2
d\mu(\tilde g(T_2)).
\eal
\leqno(A.11)
\]

Next we estimate $I_2$. From the construction, it is easy to see
that the region $\{ x \in {\bf M}, \tilde g(T_2)  \, | \, z \ge -10
\}$ is, after scaling with factor $h^{-2}$, $\delta^{1/2}$ close (in
$C^{[\delta^{-1/2}]}$ topology) to a region in the standard capped
infinite cylinder, called ${\bf C}$,  whose scalar curvature is
comparable to $1$ everywhere (see Lemma 7.3.4 [CZ] e.g.). It is well
known that the Yamabe constant of ${\bf C}$ is strictly positive,
i.e.
\[
0<Y_0 \equiv \inf_{0<v \in C^\infty({\bf C})} \frac{ 4 \frac{n-1}{n-2} \int_{\bf C} |\nabla v|^2 d\mu(g) + \int_{\bf C} R v^2
d \mu (g)}{\left( \int_{\bf C} v^{2n/(n-2)} d\mu (g) \right)^{(n-2)/n}}.
\]Here all geometric quantities are with respect to the standard metric of ${\bf C}$.

It is clear that the Yamabe constant is invariant under the scaling
of the metric. It is also known that the Yamabe constant is stable
under perturbation of metric (c.f. [BB]). Hence, when $\delta$ is
sufficiently small, we have, for  a constant $c>0$,
\[
\al
I_2&=\left( \int_{\bf M} ((1-\eta) v)^{2n/(n-2)} d\mu(\tilde g(T_2)) \right)^{(n-2)/n}\\
&\le c Y^{-1}_0
 \int_{\bf M} (
|\nabla ((1-\eta) v) |^2 + \frac{1}{4} R ((1-\eta) v)^2 ) d\mu(\tilde g(T_2)).
\eal
\leqno(A.12)
\]

The combination of (A.10), (A.11) and (A.12) gives  us
\[
\al
&\left( \int_{\bf M} v^{2n/(n-2)} d\mu(\tilde g(T_2)) \right)^{(n-2)/n}\\
&\le 2 A_1 \int (
|\nabla (\eta v) |^2 + \frac{1}{4} R (\eta v)^2 ) d\mu(\tilde g(T_2)) + 2 A_1  \int (\eta v)^2
d\mu(\tilde g(T_2))\\
&\qquad + 2c Y^{-1}_0
 \int_{\bf M} (
|\nabla ((1-\eta) v) |^2 + \frac{1}{4} R ((1-\eta) v)^2 ) d\mu(\tilde g(T_2)).
\eal
\leqno(A.13)
\]Therefore
\[
\al
&\left( \int_{\bf M} v^{2n/(n-2)} d\mu(\tilde g(T_2)) \right)^{(n-2)/n}\\
&\le 2 A_1 \int (
2 \eta^2 |\nabla v|^2 + 2 |\nabla \eta|^2 v^2 + \frac{1}{4} R (\eta v)^2 ) d\mu(\tilde g(T_2)) + 2 A_1  \int (\eta v)^2
d\mu(\tilde g(T_2))\\
&\qquad + 2c Y^{-1}_0
 \int_{\bf M} ( 2(1-\eta)^2 |\nabla v|^2 +
2 |\nabla ((1-\eta))|^2  v^2 + \frac{1}{4} R ((1-\eta) v)^2 ) d\mu(\tilde g(T_2)).
\eal
\]Combining the like terms on the right hand side of the above inequality,
we deduce
\[
\al
&\left( \int_{\bf M} v^{2n/(n-2)} d\mu(\tilde g(T_2)) \right)^{(n-2)/n}\\
&\le (4 A_1 + 4 c Y^{-1}_0) \int  (\eta^2+(1-\eta)^2) |\nabla v|^2 d\mu(\tilde g(T_2))
\\
&\qquad +
(4 A_1 + 4 c Y^{-1}_0) \int (| \nabla \eta|^2 +|\nabla (1-\eta)|^2) v^2
d\mu(\tilde g(T_2))\\
&\qquad
 + 2 A_1 \int \frac{1}{4}   \eta^2 R v^2 d\mu(\tilde g(T_2)) + 2 c Y^{-1}_0 \int
 \frac{1}{4} (1-\eta)^2 R v^2
d\mu(\tilde g(T_2)) +  2 A_1 \int v^2 d\mu(\tilde g(T_2)).
\eal
\leqno(A.14)
\]

Note that $\eta^2+(1-\eta)^2 \le 1$.  More importantly, by (A.9),
\[
| \nabla \eta|^2 +|\nabla (1-\eta)|^2 \le 2 c (R + R^-_0).
\]Also
\[
(\eta^2+(1-\eta)^2) R \le R+R^-_0.
\]Taking these inequalities to the right hand side of (A.14), we get
\[
\al
&\left( \int_{\bf M} v^{2n/(n-2)} d\mu(\tilde g(T_2)) \right)^{(n-2)/n}\\
&\le (4 A_1 + 4 c Y^{-1}_0) \int   |\nabla v|^2 d\mu(\tilde g(T_2))
\\
&\qquad +
(4 A_1 + 4 c Y^{-1}_0) 8c \int \frac{1}{4} (R+R^-_0) v^2
d\mu(\tilde g(T_2))\\
&\qquad
 + (2 A_1+ 2 c Y^{-1}_0) \int \frac{1}{4}   ( R +R^-_0) v^2 d\mu(\tilde g(T_2)) +  2 A_1 \int v^2 d\mu(\tilde g(T_2)).
\eal
\]i.e.
\[
\left( \int_{\bf M} v^{2n/(n-2)} d\mu(\tilde g(T_2)) \right)^{(n-2)/n}
\le c (A_1 + R^-_0+ Y^{-1}_0)
 \int  ( |\nabla v|^2 +\frac{1}{4} R v^2 + v^2) d\mu(\tilde g(T_2)).
 \leqno(A.15)
\]Theorem A.1 is now proven by applying Lemma A.1 and A.2. \qed

\noindent e-mail:  qizhang@math.ucr.edu
\end{document}